\documentclass{elsart}
\usepackage{amssymb}
\usepackage{amsfonts}
\usepackage{graphics}

\begin{document}

\begin{frontmatter}



\title{Graphs whose edge set can be partitioned into maximum matchings}
\begin{abstract}
This article provides structural characterization of simple graphs whose edge-set can be partitioned into maximum matchings. We use Vizing's classification of simple graphs based on edge chromatic index.
\end{abstract}
\author{Niraj Khare}
\ead{nirajkhare@math.ohio-state.edu}
\address{The Ohio State University, Columbus, Ohio, USA }
\end{frontmatter}
\section{Introduction}
By a simple graph, we shall mean a graph with no loop and no multiple edges. We will only consider simple graphs with {\bf{\underline{no isolated}}} vertex. We first fix some notations. For a graph $G$, $E(G)$ and $V(G)$ would denote the edge set and the vertex set of $G$ respectively. $\Delta(G)$, $\nu(G)$ and $\chi'(G)$ would denote the maximum degree of any vertex in $G$, the size of a maximum matching in $G$ and the edge chromatic index of $G$ respectively. For $x \in V(G)$, $deg_G(x)$ would denote the degree of the vertex $x$ and $G\setminus x$ would denote the induced subgraph on $V(G)\setminus \{x\}$.\\
We now consider simple graphs whose edge-set can be partitioned into maximum matchings. Complete graphs and even cycles are some of the examples but there are numerous other examples too.\\
Vizing's celebrated theorem states that $\chi'(G)\leq \Delta(G)+1$ for a simple graph $G$. The definition of the edge chromatic index implies that $\chi'(G)\geq \Delta(G)$. Therefore Vizing classified simple graphs as follows: a simple graph $G$ is in class I if and only if $\chi'(G)=\Delta(G)$ and a simple graph $G$  is in class II if and only if $\chi'(G)=\Delta(G)+1$. There is no structural characterization yet known for the graphs in class I or in class II.  It is NP-complete to determine whether a simple graph is in Class I or Class II (see \cite{H}). But under certain resrictions structural characterization of class I and class II graphs has been achieved. It is also known that all planar graphs with maximum degree at least seven are in class I (see \cite{SZ}). Another interesting result concerns itself with relative cardinality of class I and class II (see \cite{EW}). We will characterize Class I and class II graphs whose edge-set can be partitioned into maximum matchings.
\section{Results}
Our main aim in this article is to prove the following results.
\begin{thm} \label{thm0}
Let $G$ be a simple graph such that $\nu(G)\geq 2$, $\Delta(G)\geq 2$ and $|E(G)|=\Delta(G)\nu(G)+\lfloor{\frac{\nu(G)}{\lceil{\frac{\Delta(G)}{2}}\rceil}}\rfloor \lfloor{\frac{\Delta(G)}{2}}\rfloor$. $G$ is a unique graph up to isomorphism if and only if ${\lceil{\frac{\Delta(G)}{2}}\rceil}$ divides $\nu(G)$.
\end{thm} 
\begin{thm}\label{thm1} If $G$ is a class II graph and $E(G)$ has a partition into maximum matchings, then $\Delta(G)$ is even and $G$ is the graph with exactly $\displaystyle \frac{2\nu(G)}{\Delta(G)}$ components each isomorphic to $K_{\Delta(G)+1}$, the complete graph of order $\Delta(G)+1$.
\end{thm}
\begin{thm}\label{thm 2}
If $G$ is a class I graph and $E(G)$ can be partitioned into maximum matchings, then $E(G)$ has a partition into subgraphs that are either $K_{1,\Delta(G)}$ or a factor critical graph $H$ such that $E(H)$ can also be partitioned into maximum matchings and $\chi'(H)=\chi'(G)$.
\end{thm}
\section{Preliminaries}
We first establish some basic results that will be extremely useful in the next section. We will be borrowing some ideas and results discussed in \cite{bk}. The following definition is in \cite{bk}.
\begin{defn}\label{fdm} Let $d$, $m$ be positive integers and $G$ be a simple graph with $\Delta(G)<d$ and $\nu(G)<m$. $G \in \mathcal{F}(d,m)$ if and only if for any simple graph $G'$ such that $|E(G')|>|E(G)|$ and $G$ is a subgraph of $G'$, either $\Delta(G')\geq d$ or $\nu(G')\geq m$.
\end{defn}
\begin{thm} For a simple graph $G$, if $\nu:=\nu(G)$ and $\Delta:=\Delta(G)$, then
\begin{eqnarray}|E(G)|&\leq & \Delta\nu+\lfloor\frac{\nu}{\lceil\frac{\Delta}{2}\rceil}\rfloor \lfloor\frac{\Delta}{2}\rfloor \label{matching bar 0}\\
                      &\leq & \nu \left(\Delta+\frac{\lfloor\frac{\Delta}{2}\rfloor}{\lceil\frac{\Delta}{2}\rceil}\right).\label{matching bar 1}\end{eqnarray}.
\end{thm}
{\sl Proof.} It is obvious that for each graph $G$ with $\Delta(G)=\Delta$ and $\nu(G)=\nu$ there exists a maximal graph $G'$ such that $\nu(G')=\nu$, $\Delta(G')=\Delta$, $|E(G')|\geq |E(G)|$ and $G'\in \mathcal{F}(\Delta+1,\nu+1)$. Note that a maximal graph such as $G'$ can be constructed by adding edges or new vertices and edges to $G$. The upper bound for $|E(H)|$ obtained in \cite {bk} implies that for any $H \in \mathcal{F}(d,m)$, $|E(H)|\leq(d-1)(m-1)+\lfloor\frac{m-1}{\lceil\frac{d-1}{2}\rceil}\rfloor \lfloor\frac{d-1}{2}\rfloor$. Therefore,
\begin{eqnarray}|E(G')|&\leq & \Delta\nu+\lfloor\frac{\nu}{\lceil\frac{\Delta}{2}\rceil}\rfloor \lfloor\frac{\Delta}{2}\rfloor \\
                      &\leq & \nu \left(\Delta+\frac{\lfloor\frac{\Delta}{2}\rfloor}{\lceil\frac{\Delta}{2}\rceil}\right).\end{eqnarray}\qed\newline
Note that the above bound can also be inferred from \cite{CH}. However, the method used in \cite{CH} is more involved and it doesn't help in ascertaining uniqueness of the graphs achieving the edge bound (eq (\ref{matching bar 0})) for given $\Delta$ and $\nu$ values, i.e., Theorem \ref{thm0}. Our first goal is to prove Theorem \ref{thm0}. We next define a factor-critical graph and state the Gallai's lemma that is crucial to the following discussion.  An elegant proof of the Gallai's lemma can be found in \cite{LP}.
\begin{defn}\label{factor-critical}
A simple, connected graph $G$ is called factor-critical if and only if $G\setminus x$ has a perfect matching for all $x \in V(G)$.
\end{defn}
\begin{lem}{\bf{(Gallai)}}\label{Gallai's lemma}
Let $G$ be a simple, connected graph. If\\ $\nu(G\setminus x)=\nu(G)$ for all $x \in V(G)$, then $G$ is a factor-critical graph.
\end{lem}
We will consider graphs with \underline{no isolated vertex} only. Now we consider those graphs that attain the edge bound and analyze under what conditions these graphs are unique. We first consider some trivial cases.
\begin{rem}\label{delta1-nu1}
 For $\Delta=1$, the graph that attains the edge bound (eq (\ref{matching bar 0})) consists of $\nu$ components where each component is $K_2$. For $\nu=1$ and $\Delta >3$, the unique graph that attains the edge bound (eq (\ref{matching bar 0})) is $K_{1,\Delta}$. For $\nu=1$ and $\Delta=2$ the unique graph that attains the edge bound is $K_3$. But there are two graphs $G$ that satisfy $\Delta(G)\leq 3$, $\nu(G)=1$ and attain the edge bound $|E(G)|=3$. These graphs are $K_3$ and $K_{1,3}$. 
\end{rem}
We next consider all cases involving $\Delta\geq 2$ and $\nu \geq 2$.
\subsection{Unique graph $\mathcal{C}$}\label{graph c}
Let $\Delta \geq 2$ be a positive integer. We consider simple graphs $G$ with no isolated vertex such that
\begin{eqnarray}
\nu(G)&=&\left\lceil\frac{\Delta}{2}\right\rceil,\label{eq_3}\\
\Delta(G)&=&\Delta,\label{eq_4}\\
|E(G)|&=&\left\lfloor{\frac{(2(\lceil\frac{\Delta}{2}\rceil )+1)(\Delta)}{2}}\right\rfloor =(\Delta+\frac{\lfloor{\frac{\Delta}{2}}\rfloor}{\lceil{\frac{\Delta}{2}}\rceil})\lceil{\frac{\Delta}{2}}\rceil \label{eq_5}.\end{eqnarray}
Note that the last equation ensures that $|E(G)|$ attains the maximum edge bound given by equation (\ref{matching bar 0}). We will now construct a graph $\mathcal{C}$ satisfying equations (\ref{eq_3}), (\ref{eq_4}), (\ref{eq_5}) as follows:\\
Case (I) Let $\Delta$ be an even integer. In this case, let $\mathcal{C}=K_{\Delta+1}$.\\
Case (II) Let $\Delta =2j-1$ for some $j \geq 2$. In this case, let $\mathcal{E}$ be the graph obtained from $K_{2j}$ by removing a maximum matching. To obtain $\mathcal{C}$, we connect any $2j-1$ of the vertices of $\mathcal{E}$ to a new vertex, $v \notin V(\mathcal{E})$.\\
We next prove that $\mathcal{C}$ is the unique graph satisfying equations (\ref{eq_3}), (\ref{eq_4}) and (\ref{eq_5}).
\begin{prop}\label{unique c}
Let $\Delta \geq 2$ and $\mathcal{C}$ be the simple graph described above. If $G$ is a simple graph satisfying equations (\ref{eq_3}), (\ref{eq_4}) and (\ref{eq_5}) then\\
(a) $\nu(G\setminus x)=\nu(G)$ for all $x \in V(G)$,\\
(b) $G$ is connected,\\
(c) $G \simeq \mathcal{C}$.
\end{prop}
{\sl Proof.} Let $G$ be a graph satisfying the conditions of the proposition.\newline
{\sl Proof of (a)}: If the statement (a) is false then there exists a vertex $x \in V(G)$ such that $\nu(G\setminus x)<\nu(G)$. As at most one edge can cover $x$ in any maximum matching, we have $\nu(G\setminus x)=\nu(G)-1$. Therefore,
\begin{eqnarray*}
|E(G)| &\leq & |E(G\setminus x)|+\Delta(G\setminus x)\\
   & \leq & (\Delta(G\setminus x)+\frac{\lfloor{\frac{\Delta(G\setminus x)}{2}}\rfloor}{\lceil{\frac{\Delta(G\setminus x)}{2}}\rceil})\nu(G\setminus x)+\Delta(G\setminus x) \mbox{\quad [by equation (\ref{matching bar 1})]}\\
\end{eqnarray*}
The above expression is a non-decreasing function of $\Delta(G \setminus x)$ for a fixed $\nu(G \setminus x)$. Since $\Delta\geq \Delta(G\setminus x)$, we have\\
\begin{eqnarray*}      
|E(G)| & \leq & (\Delta+\frac{\lfloor{\frac{\Delta}{2}}\rfloor}{\lceil{\frac{\Delta}{2}}\rceil})\nu(G\setminus x)+\Delta \\
     &=& (\Delta+\frac{\lfloor{\frac{\Delta}{2}}\rfloor}{\lceil{\frac{\Delta}{2}}\rceil})(\nu-1)+\Delta\\
     &=& (\Delta+\frac{\lfloor{\frac{\Delta}{2}}\rfloor}{\lceil{\frac{\Delta}{2}}\rceil})\nu-\frac{\lfloor{\frac{\Delta}{2}}\rfloor}{\lceil{\frac{\Delta}{2}}\rceil}
\end{eqnarray*}
Therefore by assumption and the above equation, we have $|E(G)|\leq |E(G)|-\frac{\lfloor{\frac{\Delta}{2}}\rfloor}{\lceil{\frac{\Delta}{2}}\rceil}$. But $0<\frac{\lfloor{\frac{\Delta}{2}}\rfloor}{\lceil{\frac{\Delta}{2}}\rceil}$, since $\Delta \geq 2$. Hence the statement (a) holds.\newline
{\sl Proof of (b)}: On the contrary assume that $G$ is not connected. Let $\mathcal{C}_1$ be a component of $G$. Then $1\leq \nu(\mathcal{C}_1)<\nu(G)=\lceil{\frac{\Delta}{2}}\rceil $ as $G$ has no isolated vertex and $G$ is not connected by assumption. By statement (a) and Gallai's Lemma (Lemma \ref{Gallai's lemma}), $\mathcal{C}_1$ is a factor-critical component. Therefore, $|V(\mathcal{C}_1)|=2\nu(\mathcal{C}_1)+1$. So,
\begin{equation}\label{eq b}
 |E(\mathcal{C}_1)|\leq (2\nu(\mathcal{C}_1)+1)\frac{\Delta(\mathcal{C}_1)}{2}\leq (2\nu(\mathcal{C}_1)+1)\nu(\mathcal{C}_1).
\end{equation}
The above inequality implies that
\begin{eqnarray*}
\frac{|E(\mathcal{C}_1)|}{\nu(\mathcal{C}_1)} & \leq & 2 \nu(\mathcal{C}_1)+1 \\
                                              & \leq & 2(\lceil \frac{\Delta}{2} \rceil-1)+1\\
                                  &<& \Delta+ \frac{\lfloor\frac{\Delta}{2}\rfloor}{\lceil\frac{\Delta}{2}\rceil} \mbox{ \quad [as for $\Delta \geq 2$, $2 (\lceil{\frac{\Delta}{2}}\rceil-\frac{\Delta}{2})-1 \leq 0< \frac{\lfloor\frac{\Delta}{2}\rfloor}{\lceil\frac{\Delta}{2}\rceil} $].}
\end{eqnarray*}
So there is a component $\mathcal{C}_2$ of $G$ such that $\frac{|E(\mathcal{C}_2)|}{\nu(\mathcal{C}_2)}>\Delta+ \frac{\lfloor\frac{\Delta}{2}\rfloor}{\lceil\frac{\Delta}{2}\rceil}$ as $\frac{|E(G)|}{\nu(G)}=\Delta+ \frac{\lfloor\frac{\Delta}{2}\rfloor}{\lceil\frac{\Delta}{2}\rceil}$. But the equation (\ref{matching bar 1}) demands that $\frac{|E(\mathcal{C}_2)|}{\nu(\mathcal{C}_2)}\leq \Delta+ \frac{\lfloor\frac{\Delta}{2}\rfloor}{\lceil\frac{\Delta}{2}\rceil}$. The contradiction implies that the statement (b) holds.\newline
{\sl Proof of (c)}: Since the statements (a) and (b) hold for $G$, $G$ is factor-critical by Gallai's Lemma (Lemma \ref{Gallai's lemma}). As $\nu(G)=\lceil{\frac{\Delta}{2}}\rceil$, we have $|V(G)|=2(\lceil{\frac{\Delta}{2}}\rceil)+1$. We consider following two cases.\newline
If $\Delta$ is even then $G$ is a connected graph with $2(\lceil{\frac{\Delta}{2}}\rceil)+1=\Delta+1$ vertices and $|E(G)|=(\Delta+ \frac{\lfloor\frac{\Delta}{2}\rfloor}{\lceil\frac{\Delta}{2}\rceil})\lceil\frac{\Delta}{2}\rceil= \frac{(\Delta+1)\Delta}{2}$. Therefore, $deg_G(x)=\Delta$ for all $x \in V(G)$. Hence $G$ is $K_{\Delta+1}$, the complete graph on $\Delta+1$ vertices. So $G\simeq \mathcal{C}$.\newline
If $\Delta$ is odd. Let $\Delta=2j-1$ for some $j\geq 2$. Then
 $\nu({G})=\left\lceil\frac{\Delta}{2}\right\rceil=j$, 
$|V(G)|=2\nu(G)+1=2j+1$ and 
$|E({G})|=(\Delta+ \frac{\lfloor\frac{\Delta}{2}\rfloor}{\lceil\frac{\Delta}{2}\rceil})\lceil\frac{\Delta}{2}\rceil=(2j-1)j+j-1.$
So $$\sum_{x \in V(G)}deg_G(x)=2j(2j-1)+2j-2.$$ Therefore there is a unique vertex $v \in V({G})$ of degree $2j-2$. Hence there is a vertex $u$ in $V({G})$ which is not a neighbor of $v$. Consequently ${G} \setminus u$ is a regular graph of degree $2j-2$ on $2j$ vertices and hence its complement is a regular graph of degree one, namely, a matching of a complete graph on $2j$ vertices. This establishes ${G}\simeq \mathcal{C}$, where $\mathcal{C}$ is the graph described earlier for the case $\Delta=2j-1$.\newline
\subsection{Unique graphs with maximum number of edges for a given maximum degree and matching size}
We emphasize that graphs considered in this discussion have no isolated vertex.  Note that a method is provided in \cite{bk} to construct a graph $G$ such that $|E(G)|$ attains the edge bound given by equation (\ref{matching bar 0}). We find the condition when the graphs that attain the maximum edge bound (equation (\ref{matching bar 0})) are unique up to isomorphism.
\begin{prop}\label{unique_graphs} Let $\mathcal{C}$ be the graph constructed in the subsection \ref{graph c} and $G$ be a simple graph with $\nu:=\nu(G)$ and $\Delta:=\Delta(G)$ such that $\Delta\geq 2$. If $\lceil{\frac{\Delta}{2}}\rceil$ divides $\nu$ and $\displaystyle |E(G)|=\Delta\nu+\frac{\nu\left\lfloor\frac{\Delta}{2}\right\rfloor}{\lceil\frac{\Delta}{2}\rceil}$, then \\
(a) $\nu(G\setminus x)=\nu(G)$ for all $x \in V(G)$,\\
(b) if $D$ is a component of $G$ with $\nu(D)> \nu(\mathcal{C})$, then $\frac{|E(D)|}{\nu(D)} < \frac{|E(\mathcal{C})|}{\nu(\mathcal{C})}$,\\
(c) if $D$ is a component of $G$ with $\nu(D) <\nu(\mathcal{C})$, then $\frac{|E(D)|}{\nu(D)} < \frac{|E(\mathcal{C})|}{\nu(\mathcal{C})} $,\\
(d) every component of $G$ is isomorphic to $\mathcal{C}$,\\
(e) $G$ is unique up to isomorphism.
\end{prop}
{\sl Proof.} Let $G$ be a graph satisfying the conditions of the proposition.\newline
{\sl Proof of (a)}: If the statement (a) is false then there exists a vertex $v \in V(G)$ such that $\nu(G\setminus v)<\nu(G)$. This implies
\begin{eqnarray*}
\Delta\nu+\frac{\nu}{\lceil{\frac{\Delta}{2}}\rceil} \lfloor\frac{\Delta}{2}\rfloor=|E(G)|&\leq& deg_G(v)+|E(G\setminus v)|\\
                                  &\leq& \Delta+ \Delta(G\setminus v)(\nu-1)+\left\lfloor\frac{\nu-1}{\lceil\frac{\Delta(G\setminus v)}{2}\rceil}\right\rfloor \left\lfloor\frac{\Delta(G\setminus v)}{2}\right\rfloor
\end{eqnarray*}
We again note that the edge bound, i.e., the equation (\ref{matching bar 0}) is a non-decreasing function of $\Delta$ for a fixed $\nu$. Also $\Delta\geq \Delta(G \setminus v)$. Therefore,
\begin{eqnarray*}
\Rightarrow\quad \frac{\nu}{\lceil\frac{\Delta}{2}\rceil}{\left\lfloor \frac{\Delta}{2}\right\rfloor}\leq \left\lfloor \frac{\nu-1}{\lceil\frac{\Delta}{2}\rceil}\right\rfloor{\left\lfloor \frac{\Delta}{2}\right\rfloor} &\leq& \frac{(\nu-1)}{\lceil\frac{\Delta}{2}\rceil}{\left\lfloor \frac{\Delta}{2}\right\rfloor}\\
\Rightarrow\qquad \nu &\leq& \nu-1.
\end{eqnarray*}
This contradiction proves (a).\newline
We recall that the graph $\mathcal{C}$ is a factor critical graph with $\Delta \lceil\frac{\Delta}{2}\rceil+\lfloor\frac{\Delta}{2}\rfloor$ edges and maximum matching size $\lceil\frac{\Delta}{2}\rceil$.\newline
{\sl Proof of (b)}: Let $D$ be a component of $G$. Gallai's lemma and (a) imply that $D$ is factor-critical, hence $|V(D)|=2\nu(D)+1$ and
\[|E(D)| \leq \left\lfloor\frac{(2\nu(D)+1)\Delta}{2}\right\rfloor = \nu(D)\Delta +\left\lfloor\frac{\Delta}{2}\right\rfloor.\]
Note that in the following inequality we again used the fact that the edge bound, i.e., the equation (\ref{matching bar 0}) is a non-decreasing function of $\Delta$ for a fixed $\nu$. If $\nu(D)> \nu(\mathcal{C})$, then \[\quad \frac{|E(D)|}{\nu(D)} \leq \Delta +\frac{\lfloor\frac{\Delta}{2}\rfloor}{\nu(D)}< \Delta +\frac{\lfloor\frac{\Delta}{2}\rfloor}{\nu(\mathcal{C})}=\Delta+\frac{\lfloor\frac{\Delta}{2}\rfloor}{\lceil\frac{\Delta}{2}\rceil}=\frac{|E(\mathcal{C})|}{\nu(\mathcal{C})}.\]
This proves (b).\newline
{\sl Proof of (c)}: Let $D$ be a component of $G$. 
By Gallai's lemma and (a), $D$ is a factor critical graph and hence $|V(D)|=2\nu(D)+1$. Thus 
\[|E(D)|\leq \frac{(2\nu(D)+1)2\nu(D)}{2}=(2\nu(D)+1)\nu(D).\]
Now since $\nu(D) \le \nu(\mathcal{C})-1$, we have 
\[\frac{|E(D)|}{\nu(D)} \le 2\nu(D)+1 \le 2\nu(\mathcal{C})-1=2\left\lceil\frac{\Delta}{2}\right\rceil-1<\Delta+\frac{\lfloor\frac{\Delta}{2}\rfloor}{\lceil\frac{\Delta}{2}\rceil}=\frac{|E(\mathcal{C})}{\nu(\mathcal{C})},\]
which proves (c).\newline
{\sl Proof of (d)}: From (b) and (c), every component $D$ of $G$ has $\nu(D)=\nu(\mathcal{C})=\lceil\frac{\Delta}{2}\rceil$ and also $|E(D)|=\Delta \lceil\frac{\Delta}{2}\rceil+\lfloor\frac{\Delta}{2}\rfloor$. Since $\mathcal{C}$ is the only graph with maximum matching size $\lceil\frac{\Delta}{2}\rceil$ and the number of edges $\Delta \lceil\frac{\Delta}{2}\rceil+\lfloor\frac{\Delta}{2}\rfloor$, $D$ must be isomorphic to $\mathcal{C}$. This proves (d).\newline
{\sl Proof of (e)}: This follows from (d).\qed

Now we explore the inverse of Proposition \ref{unique_graphs}. Recall the connected graph $\mathcal{C}$ described in the subsection \ref{graph c} is the unique graph satisfying equations (\ref{eq_3}), (\ref{eq_4}) and (\ref{eq_5}).
\begin{prop}\label{inverse_unique_graphs}
Let $G_1$ be a simple graph such that $|E(G_1)|=\Delta(G_1)\nu(G_1)+\lfloor{\frac{\nu(G_1)}{\lceil{\frac{\Delta(G_1)}{2}}\rceil}}\rfloor \lfloor{\frac{\Delta(G_1)}{2}}\rfloor $, i.e., $G_1$ attains the maximum edge bound given by the inequality (\ref{matching bar 0}). If $\nu(G_1)\geq 2$, $\Delta(G_1)\geq 2$ and ${\lceil{\frac{\Delta(G_1)}{2}}\rceil}$ doesn't divide $\nu(G_1)$, then there exists a simple graph $G_2$ such that $|E(G_2)|=|E(G_1)|$, $\nu(G_2)=\nu(G_1)$, $\Delta(G_2)=\Delta(G_1)$ and $G_2$ is \underline{not} isomorphic to $G_1$.  
\end{prop}
{\sl Proof.} We use the method given in \cite{bk} to construct a simple graph $G$ such that $|E(G)|=|E(G_1)|$, $\nu(G)=\nu(G_1)$ and $\Delta(G)=\Delta(G_1)$. Let $G$ have $t:=\nu(G_1)-({\lceil{\frac{\Delta(G_1)}{2}}\rceil})\lfloor{\frac{\nu(G_1)}{{\lceil{\frac{\Delta(G_1)}{2}}\rceil}}}\rfloor$ components isomorphic to $K_{1,\Delta(G_1)}$ and $\lfloor{\frac{\nu(G_1)}{{\lceil{\frac{\Delta(G_1)}{2}}\rceil}}}\rfloor$ components isomorphic to $\mathcal{C}$ (described in subsection \ref{graph c}). Let $G_2=G$ if $G_1$ is \underline{not} isomorphic to $G$. So assume that $G_1$ is isomorphic to $G$. Note that $t:=\nu(G_1)-({\lceil{\frac{\Delta(G_1)}{2}}\rceil})\lfloor{\frac{\nu(G_1)}{{\lceil{\frac{\Delta(G_1)}{2}}\rceil}}}\rfloor\geq 1$ as ${\lceil{\frac{\Delta(G_1)}{2}}\rceil}$ doesn't divide $\nu(G_1)$. If $t\geq 2$ then let $H_1$ and $H_2$ be two components of $G_1$ isomorphic to $K_{1,\Delta(G_1)}$. We remove an edge of $H_1$ and then connect the vertex of degree $\Delta(G_1)-1$ of $H_1$ to any vertex of degree one in $H_2$ by a new edge. Thus, a new graph $G_2$ is obtained. It is obvious by construction that $G_2$ satisfies hypothesis of Proposition \ref{inverse_unique_graphs}. So, we need to consider only the case $t=1$ to complete the proof. As $\nu(G_1)\geq 2$ and $t=1$, $G_1$ has at least two components. Hence $G_1$ has a component isomorphic to $K_{1,\Delta(G_1)}$ and at least one component isomorphic to $\mathcal{C}$. As $\Delta(G_1)+|E(\mathcal{C})|=\Delta(G_1)+\Delta(G_1)(\lceil{\frac{\Delta(G_1)}{2}}\rceil)+\lfloor{\frac{\Delta(G_1)}{2}}\rfloor \leq \lfloor{\frac{(2(\lceil\frac{\Delta(G_1)}{2}\rceil +1)+1)\Delta(G_1)}{2}}\rfloor$, we can coalesce the two components to form a factor critical component with $2(\lceil\frac{\Delta(G_1)}{2}\rceil +1)+1$ vertices and maximum degree at most $\Delta(G_1)$ which has number of edges equal to $\Delta(G_1)+|E(\mathcal{C})|$ and maximum matching size equal to $\nu(\mathcal{C})+1$. Thus, again a new graph $G_2$ is obtained. It is obvious by construction that $G_2$ satisfies hypothesis of Proposition \ref{inverse_unique_graphs}.\qed

We can combine the above two propositions in the following theorem. 
\begin{thm}\label{uniqueness_thm} 
Let $G$ be a simple graph such that $\nu(G)\geq 2$, $\Delta(G)\geq 2$ and $|E(G)|=\Delta(G)\nu(G)+\lfloor{\frac{\nu(G)}{\lceil{\frac{\Delta(G)}{2}}\rceil}}\rfloor \lfloor{\frac{\Delta(G)}{2}}\rfloor$. $G$ is a unique graph up to isomorphism if and only if ${\lceil{\frac{\Delta(G)}{2}}\rceil}$ divides $\nu(G)$.
\end{thm}
{\sl Proof.} The above conclusion follows by Proposition \ref{unique_graphs} and Proposition \ref{inverse_unique_graphs}.\qed
\section{Graphs whose edge set can be partitioned into maximum matchings} 
\begin{defn}\label{friendly_colorable} A simple graph $G$ is called \textit{friendly-edge-colorable} if and only if $E(G)$ has a partition into maximum matchings.
\end{defn}
Examples: $C_4$, $C_{2m}$ for $m\geq 2$, $K_{2m}$ for $m\geq 1$, $K_{2m+1}$ for $m\geq 1$ and $K_{1,n}$ for $n\geq 1$.\newline
We observe that a proper edge coloring of $G$ is equivalent to partitioning $E(G)$ into matchings (may not be maximum).
\begin{prop}\label{f_ob_1}
$G$ is a friendly-edge-colorable graph if and only if $|E(G)|=\chi'(G)\nu(G)$.
\end{prop}
{\sl Proof.} We consider a minimal proper edge coloring of $G$. Since every color class of $G$ is a matching and a matching in $G$ can be of size at most $\nu(G)$, we have $|E(G)|\leq \chi'(G) \nu(G)$.
So if $G$ is a friendly-edge-colorable graph, then there is a partition of $E(G)$ into maximum matchings. Hence there is a positive integer $n$ such that $|E(G)|=n\nu(G)$. This partition corresponds to a proper edge coloring with $n$ colors. Thus, $\chi'(G)\leq n$. Hence,
$|E(G)|=n \nu(G)\ge \chi'(G)\nu(G)$. Thus for a friendly-edge-colorable graph $G$, we have  $|E(G)| = \chi'(G)\nu(G)$.\newline
Suppose $G$ is \underline{not} a friendly-edge-colorable graph. We consider an edge coloring of $G$ in $\chi'(G)$ colors. Note that each color class is a matching and at least one of the color classes is of size strictly less than $\nu(G)$  as $G$ is not a friendly-edge-colorable graph. Thus, we get $|E(G)|<\chi'(G)\nu(G)$.\qed 
\begin{prop}\label{f_ob_2}
If $G$ is a friendly-edge-colorable graph, then any proper edge coloring of $G$ in $\chi'(G)$ colors results in a partition of $E(G)$ into maximum matchings.
\end{prop}
\noindent{\sl Proof.} Any proper edge coloring of $G$ is a partition of $E(G)$ into matchings (may not be maximum). Consider a proper edge coloring of $G$ in $\chi'(G)$ colors. If there is a color class of size strictly less than $\nu$, then $|E(G)|<\chi'(G)\nu(G)$ which contradicts Proposition \ref{f_ob_1}.\qed
\begin{rem} The name ``friendly-edge-colorable" is due to Proposition \ref{f_ob_2} which states that when the least number of colors are used to properly color the edges, the colors are equally distributed so that each color class gets the maximum!
\end{rem}
The following theorem characterizes friendly-edge-colorable graphs in class II.
\begin{thm}\label{f_thm_1} Let $G$ be a friendly-edge-colorable graph of class II such that $\Delta(G)\geq 2$ and $\nu(G)\geq 2$. Then\\
(a) $\Delta(G)$ is even,\\
(b) $\displaystyle \frac{\Delta(G)}{2}$ divides $\nu(G)$,\\
(c) every component of $G$ is isomorphic to $K_{\Delta(G)+1}$, the complete graph of order $\Delta(G)+1$.
\end{thm}
{\sl Proof.} Let $\Delta=\Delta(G)$ and $\nu=\nu(G)$. By Proposition \ref{f_ob_2}, $|E(G)|=\nu(\Delta+1)$. Also by inequality (\ref{matching bar 0}),
\begin{eqnarray*}
(\Delta+1)\nu=|E(G)|&\leq & \Delta\nu+\lfloor\frac{\nu}{\lceil\frac{\Delta}{2}\rceil}\rfloor \lfloor\frac{\Delta}{2}\rfloor,\\
                     &\leq& \nu \Delta + \frac{\nu}{\lceil\frac{\Delta}{2}\rceil} \lfloor\frac{\Delta}{2}\rfloor\\
                     & \leq & \nu \Delta+ \nu, \mbox{[ as $\frac{\lfloor\frac{\Delta}{2}\rfloor}{\lceil\frac{\Delta}{2}\rceil} \leq 1$].}\\
\end{eqnarray*}
Thus, $\displaystyle \left\lfloor\frac{\nu}{\lceil\frac{\Delta}{2}\rceil}\right\rfloor=\frac{\nu}{\lceil\frac{\Delta}{2}\rceil}$, which proves (b). Also, $\displaystyle \frac{\lfloor{\frac{\Delta}{2}}\rfloor}{\lceil{\frac{\Delta}{2}}\rceil}=1$ that implies (a). By Proposition \ref{unique_graphs}, we know that every component of $G$ is isomorphic to $K_{\Delta+1}$, the complete graph of order $\Delta+1$.\qed 

Now we will consider friendly-edge-colorable graphs that are in class I, i.e., graphs with edge chromatic index $\Delta$. Next we prove three lemmas that help us characterize friendly-edge-colorable graphs in class I.
\begin{lem}\label{f_lem_1}
Let $G$ be a friendly-edge-colorable graph. If there exists a vertex $x \in V(G)$ such that $\nu(G\setminus x)<\nu(G)$ then $deg_G(x)=\chi'(G)$.
\end{lem}
\noindent{\sl Proof.} As $G$ is friendly-edge-colorable, there exists a partition of $E(G)$ into maximum matchings and, by Proposition \ref{f_ob_1}, $\frac{|E(G)|}{\nu(G)}=\chi'(G)$. Hence each class, i.e., each part of this partition, has size $\nu(G)$ and there are $\chi'(G)$ parts. We denote each color class (i.e., a part in the partion) by $C_i$ for all $i\in \{1, \ldots, \chi'(G)\}$. If $deg_G(x)<\chi'(G)$ then there is a color missing at $x$. Without loss of generality, let the missing color be $1$. Hence the whole color class $C_1$ belongs to $E(G\setminus x)$, i.e., $C_1\subseteq E(G\setminus x)$. Therefore, $\nu(G\setminus x)\geq |C_1|$. But $C_1$ is a maximum matching of $G$ hence $|C_1|=\nu(G)$. This contradicts that $\nu(G\setminus x)<\nu(G)$. Hence $deg_G(x)\geq \chi'(G)$. Also $deg_G(x)\leq \Delta(G)\leq \chi'(G)$ by the definition of $\chi'(G)$.\hspace{\stretch{1}} $\square$\newline   
\begin{rem}\label{f_rem_1}
Lemma \ref{f_lem_1} implies that if $G$ is friendly-edge-colorable and there exists $x \in V(G)$ such that $\nu(G\setminus x)<\nu(G)$, then $deg_G(x)=\Delta(G)=\chi'(G)$ as $\Delta(G)\leq\chi'(G)$. Hence $G$ is in class I. So the contrapositive implies that if $G$ is friendly-edge-colorable and in class II, i.e. \underline{not} in class I, then $\nu(G\setminus x)=\nu(G)$ for all $ x \in V(G)$. By Gallai's lemma, each component of $G$ is a factor-critical component as we noticed in Theorem \ref{f_thm_1}.
\end{rem}
\begin{lem}\label{f_lem_2}
Let $G$ be a friendly-edge-colorable graph. If there exists $ x \in V(G)$ such that $\nu(G\setminus x)<\nu(G)$ and $|E(G\setminus x)|\geq 1$, then $G\setminus x$ is a friendly-edge-colorable graph and $\chi'(G\setminus x)=\chi'(G)$.
\end{lem}
\noindent{\sl Proof.} As $G$ is friendly-edge-colorable, $E(G)$ has a partition into maximum matchings. Define a coloring, $C$, corresponding to such a partition. Then each color class has size $\nu(G)$ and there are $\frac{|E(G)|}{\nu(G)}=\chi'(G)$ color classes by Proposition \ref{f_ob_1}. Now consider the restriction of $C$ on $E(G\setminus x)$. By Lemma \ref{f_lem_1}, $deg_G(x)=\chi'(G)$ hence $C_{E(G\setminus x)}$ has $\chi'(G)$ color classes each of size $\nu(G)-1$. As by assumption $\nu(G\setminus x)=\nu(G)-1$, $E(G\setminus x)$ has a partition into maximum matchings, implying that $G\setminus x$ is friendly-edge-colorable. Since $G\setminus x$ is friendly-edge-colorable, by Proposition \ref{f_ob_1},we get
\begin{eqnarray*}
\chi'(G\setminus x) &=& \frac{|E(G\setminus x)|}{\nu(G\setminus x)}\\
                    &=& \frac{|E(G)|-deg_G(x)}{\nu(G)-1}\\   
                    &=& \frac{|E(G)|-\chi'(G)}{\nu(G)-1}\\
                    &=& \frac{\nu(G)\chi'(G)-\chi'(G)}{\nu(G)-1}\\
                    &=& \chi'(G).
\end{eqnarray*}\hspace{\stretch{1}} $\square$\newline
By a non-trivial component of a simple graph, we shall mean a component that has at least an edge, i.e., $K_2$ is a subgraph of the component.
\begin{lem}\label{f_lem_3}
$G$ is a friendly-edge-colorable graph if and only if each non-trivial component $C$ of $G$ is friendly-edge-colorable and $\chi'(G)=\chi'(C)$.
\end{lem}
\noindent{\sl Proof.} We first show the if part. Consider a proper edge coloring for each component in $\chi'(G)$ colors. This gives the desired partition of the edge set of $G$.\newline
Next we show the only if part. Since $G$ is friendly-edge-colorable, there is a partition of $E(G)$ into maximum matchings that has $\chi'(G)$ classes (by Proposition \ref{f_ob_1}). We claim that restricting this partition to any component of $G$ gives a partition of the edge set of the component into maximum matchings of the component. Suppose on the contrary that some matching $\mathcal{M}_{\mathcal{C}}$ obtained this way (by restricting a maximum matching $\mathcal{M}$ to $\mathcal{C}$) for a component $\mathcal{C}$ is not a maximum matching in $\mathcal{C}$. Then there is an augmenting path relative to this matching in the component. But then this is an augmenting path, relative to the corresponding matching $\mathcal{M}$ of $G$, in $G$. Hence $\mathcal{M}$ is not a maximum matching of $G$. This contradicts the fact that $\mathcal{M}$ is a maximum matching.\newline
If for some component $C$, $\chi'(C)\neq \chi'(G)$ then $\chi'(C)<\chi'(G)$ (as $\chi'(C)\leq\chi'(G)$ for $C$ is a component of $G$). Now consider any proper edge coloring of $G$ in $\chi'(G)$ colors and for each component order color classes greedily. For each $1\leq i \leq \chi'(G)$, we denote the $i-{th}$ color class by $\mathcal{D}(i)$. Then the first $\chi'(\mathcal{C})$ color classes will have strictly larger size than the remaining color classes $\{ \mathcal{D}(\chi'(\mathcal{C})+1), \mathcal{D}(\chi'(\mathcal{C})+2),\ldots,\mathcal{D}(\chi'(G))\}$. This contradicts that each color class is of the same size $\nu(G)$ as shown by Proposition \ref{f_ob_2}.\hspace{\stretch{1}} $\square$
\begin{rem}\label{f_rem_2}
In the previous lemma the condition that $\chi'(G)=\chi'(C)$ for all components $C$ is necessary otherwise counter examples exist. For instance, let $G$ be the graph consisting of two components one of which is $K_3$ and the other one is $K_{1,4}$.
\end{rem}
\begin{thm}\label{f_thm_2}
If $G$ is a non-trivial, class I, friendly-edge-colorable graph then $E(G)$ has a partition into the following two kinds of subgraphs:\newline
(i) $K_{1,\Delta(G)}$\newline
(ii) factor critical, friendly-edge-colorable graphs with edge chromatic index $\chi'(G)$.
\end{thm}
\noindent{\sl Proof.} Let $G_0:=G$. If there exists $x_1 \in V(G)$ such that $\nu(G\setminus x_1)<\nu(G)$ then remove all edges incident to $x_1$ and define $G_1:=G\setminus x_1$. Since by Lemma \ref{f_lem_1} and Remark \ref{f_rem_1} $deg_G(x_1)=\Delta$, we removed a $K_{1,\Delta(G)}$ from $G$. Note that $G_1$ is a friendly-edge -colorable graph and $\chi'(G_1)=\chi'(G)$ by Lemma \ref{f_lem_2}. By Lemma \ref{f_lem_3}, each component of $G_1$ is also a friendly-edge-colorable graph and has edge chromatic index $\chi'(G)$. Similarly, for $i\geq 1$ define $G_i:=G_{i-1}\setminus x_i$ if there exists $x_i \in V(G_{i-1})$ such that $\nu(G_{i-1}\setminus x_i)<\nu(G_{i-1})$. We remove the vertices $x_i$ and all edges adjacent to each of these vertices. Corresponding to each of the $x_i$'s, we get a subgraph isomorphic to $K_{1,\Delta(G)}$. For some large enough $i$, we have only those non-trivial components $C$ in $G_i$ such that $\nu(C\setminus x)=\nu(C)$ for all $ x\in V(C)$. By Gallai's lemma, these are the factor critical components. Also by Lemma \ref{f_lem_3}, $\chi'(C)=\chi'(G_i)=\chi'(G)$.\hspace{\stretch{1}} $\square$
\begin{rem} Reader can review Figure \ref{fig:max_match_3} and notice that the degree of the vertex $x$ in Figure \ref{fig:max_match_3} is $3$ and $x$ must be covered by every maximum matching. Also, removal of the vertex $x$ from the graph yields a graph whose only non-trivial component is factor-critical.
\end{rem}
\section{Acknowledgements}
I would like to thank Dr. Nishali Mehta, Dr. Naushad Puliyambalath and Prof. \'Akos Seress for their valuable comments and help to improve the draft.    

\end{document}